\documentclass[11pt]{article}
\usepackage{amsfonts,amsbsy,amssymb,amsmath,graphicx}
\usepackage{epsfig}
\usepackage{dsfont}
\usepackage{mathrsfs}
\usepackage{geometry}
\usepackage{float}
\usepackage{multirow}
 \usepackage[autostyle=true]{csquotes}
\usepackage{booktabs}
\usepackage{subcaption}
\usepackage{enumitem}
\usepackage{color}
\usepackage{multirow}
\usepackage{color}
\usepackage{natbib}
\usepackage{rotating}
\usepackage{epsfig}
\usepackage{ifpdf}
\usepackage[colorlinks=true,linkcolor=red,citecolor=blue]{hyperref}

\def\1{1\kern-.20em {\rm l}}

\newtheorem{theorem}{Theorem}

\newtheorem{lemma}{Lemma}

\newcommand{\R}{{\mathbb{R}}}

\newcommand{\Pb}{{\mathbb{P}}}
\newcommand{\Eb}{{\mathbb{E}}}

\newcommand{\var}{\operatorname{Var}}

\newcommand{\ft}{^{ \text{ft} } }

\DeclareMathOperator{\sign}{sign}

\DeclareMathOperator{\const}{const}

\newcommand{\rw}{\widehat{r}}
\newcommand{\fw}{\widehat{f}}
\newcommand{\mw}{\widehat{m}}

\numberwithin{equation}{section}

\def\1{1\kern-.20em {\rm l}}

\numberwithin{equation}{section}

\begin{document}


\title{\bf  Partial heteroscedastic deconvolution estimation in nonparametric regression}

\author{ {\sc Baba Thiam}\footnote{email: baba.thiam@univ-lille.fr} \\
Univ. Lille, CNRS, UMR 8524 - Laboratoire Paul Painlevé, F-59000 Lille, France\\
}

\maketitle

\begin{abstract}
\noindent In this paper, we consider a partial deconvolution kernel estimator for nonparametric regression when some covariates are measured with error while others are observed without error. We focus on a general and realistic setting in which the measurement errors are heteroscedastic. We propose a kernel-based estimator of the regression function in this framework and show that it achieves the optimal convergence rate under suitable regularity conditions. The finite-sample performance of the proposed estimator is illustrated through simulation studies.
\noindent   

\vspace{1mm}

\end{abstract}

\section{Introduction}\label{intro}
In statistical analysis, it is common to investigate the relationship between explanatory variables $X$ and a response variable $Y$ using regression techniques. Regression estimation can be carried out nonparametrically through smoothing methods such as the Nadaraya–Watson estimator, local polynomial smoothing, or smoothing splines, among others. These methods are particularly effective when all explanatory variables are directly observable. However, in many practical situations, some variables cannot be observed precisely and are instead measured with error. Such measurement errors may arise from imperfect measuring instruments or from the inherent difficulty of accessing the variable of interest. Measurement error occurs in numerous applied fields, including astronomy, nutrition, epidemiology, chemistry, and many others. When all explanatory variables are contaminated by measurement error, a wide range of methodologies has been developed to reduce its impact, particularly in parametric regression models under specific distributional assumptions on the error. Comprehensive studies of these issues can be found in \cite{F87} in the framework of linear regression models while  \cite{CRS12} provide an extensive overview for nonlinear regression models. The latter also provide an overview of some methodological developments of significant applied relevance, focusing in particular on deconvolution methods for the nonparametric estimation of density and regression functions affected by measurement error. When all covariates are measured with error, the estimation of probability density functions and regression functions has received considerable attention in the literature. Large-sample properties of deconvolution kernel density estimators have been studied, for example, by \cite{SC90}, \cite{CH88}, \cite{LT89}, \cite{F91a, F91b, F91c}. The asymptotic properties of deconvolution kernel regression estimators were established by \cite{FT93}who showed that the optimal local and global convergence rates depend on the tail behavior of the characteristic function of the measurement error distribution. As emphasized by \cite{SBS20}, in many real-data applications some explanatory variables may be observed without error, while others are contaminated by measurement error. A classical example arises in the Framingham Heart Study, discussed in \cite{Kan86}, where the goal is to understand how the development of coronary heart disease depends on factors such as systolic blood pressure, age, body mass, serum cholesterol, and smoking status. In this study, age, body mass, and smoking status are measured accurately, whereas serum cholesterol and systolic blood pressure are subject to measurement error. In such settings, it is natural to replace standard kernel functions with deconvolution kernels for the covariates affected by measurement error. More formally, let $X$ denote the observable explanatory variables and $T$ the unobservable covariates measured with error. Let $W$ be the surrogate variable defined by $W=T+U$. We consider the regression model
\begin{eqnarray}
\label{model1}
Y=r(X, T)+\varepsilon, ~~~ W=T+U, ~~~\mathbb{E}(\varepsilon|X,T)=0,
\end{eqnarray}
where $X$ and $T$ have densities $f_X$ and $f_T$, respectively,  and the measurement error $U$ is independent of $(X,T,\varepsilon)$ with know  density $f_U$, where the error variable $U$ isand the error probability density function $f_U$ is known. For simplicity we assume that $X$ and $T$ are independent and let  $f_{(X,T)}$ be the joint density of $(X,T)$.\\
Suppose that we observe a sample of independent and identically distributed (iid) random vectors $(X_1,T_1,Y_1), \ldots, (X_n,T_n,Y_n)$ generated from model \eqref{model1}. To estimate the regression function $r$, \cite{SBS20} proposed the following estimator
\begin{eqnarray}
\label{estim1}
\rw_n(x,t)=(nhb)^{-1}\sum_{j=1}^{n}K\left(\frac{x-X_j}{h}\right)L_U\left(\frac{t-W_j}{b}\right)/\widehat{f}_n(x,t),
\end{eqnarray}
where 
$\fw_n(x,t)=(nhb)^{-1}\sum_{j=1}^{n}Y_{j}K\left(\frac{x-X_j}{h}\right)L_U\left(\frac{t-W_j}{b}\right)$ is the partial deconvolution kernel density estimator of $f_{(X,T)}$, $K$ and $L$ are kernels functions,
$$L_U(t)=(2\pi)^{-1}\textstyle\int e^{-itv}L^{\text{ft}}(v)/f_U^{\text{ft}}(v/b)dv,~ \text{with}~i^2=-1,$$and
$g^{\text{ft}}$ denotes the Fourier transform of a function $g$.
In \cite{SBS20}, the authors established the optimal local and global convergence rates for the estimator $\widehat{r}_n(x,t)$ under both ordinary smooth and super smooth error distributions. \\
However, as emphasized in \cite{DM07}, the assumption of homoscedastic measurement errors $U_j$ is often unrealistic in practical applications. In many real-data settings, heteroscedasticity in measurement errors arises when data are collected under heterogeneous conditions. For instance, datasets may be constructed by aggregating measurements from multiple laboratories or research centers (see, e.g., \cite{NRC93}), or by combining results from distinct studies, as commonly encountered in meta-analytic frameworks (see, e.g., \cite{W97}). Moreover, different subpopulations (such as healthy versus unhealthy individuals, or smokers versus non-smokers) may be subject to distinct contamination mechanisms, as already discussed by \cite{F87}. Finally, the measurement process itself may be subjective and vary across individuals; for example, \cite{BF54} documented substantial discrepancies in students’ subjective assessments of iron content in various substances.\\
Motivated by these considerations, we consider the following heteroscedastic error-in-variables model:
\begin{eqnarray}
\label{mod_het}
Y_j=r(X_j,T_j)+\varepsilon_j, ~~~W_j=T_j+U_j, ~~~\mathbb{E}(\varepsilon_j|X_j,T_j)=0,
\end{eqnarray}
with $X_j\leadsto f_X$, $T_j\leadsto f_T$ and $U_j\leadsto f_{U_j}$.  The measurement error $U_j$ is assumed to be independent of $(X_j,T_j, Y_j,\varepsilon_j)$. Furthermore, we assume that $X_j$ and $T_j$ are independent. The error densities $f_{U_j}$ may depend both on the observation index $j$ and on the sample size $n$. 
As pointed out in \cite{DM07}, the estimator \eqref{estim1} cannot be directly applied in this setting, since it relies on a single (common) error density. Therefore, considering observations generated from the heteroscedastic error model \eqref{mod_het} and following the approach of \cite{SBS20}, we propose the following partial deconvolution regression estimator:
\begin{eqnarray}
\label{estim}
\widehat{r}_n(x,t)=\frac{1}{hb}\sum_{j=1}^{n}Y_{j}K\left(\frac{x-X_j}{h}\right)L_{U_j}\left(\frac{t-W_j}{b}\right)/\widehat{f}_n(x,t),
\end{eqnarray}
where 
\begin{eqnarray}
\label{denom}
\widehat{f}_n(x,t)&=&\frac{1}{hb}\sum_{j=1}^{n}K\left(\frac{x-X_j}{h}\right)L_{U_j}\left(\frac{t-W_j}{b}\right)
\end{eqnarray}
is an estimator of the joint density $f_{(X,T)}$ that remains valid under heteroscedastic measurement errors, with
\begin{eqnarray}
\label{lu_j}
L_{U_j}(t)=\frac{1}{2\pi}\int_{\R} \exp(-ivt)L\ft(v)\psi_j(v/b)dv,
\end{eqnarray}
and where
\begin{eqnarray}
\label{psi_j}
\psi_j(v)=f_{U_j}^{\text{th}}(-v)/(\sum_{k=1}^n|f_{U_k}^{\text{th}}(v)|^2).
\end{eqnarray}
Note that the function $\psi_j$ generalizes the factor $\left(nf_U\ft\right)^{-1}$ used in the homoscedastic framework considered in \cite{SBS20}. \\
The  remainder of the paper is organized as follows. Section \ref{ass-res} presents the assumptions and theoretical properties of the proposed estimator. Section~\ref{numeric} reports numerical results illustrating its finite-sample performance. Finally, Section~\ref{proofs} is devoted to the technical proofs.

\section{Assumptions and main results}\label{ass-res}
Define $\tau^2(x,t)=\var\left(Y|X=x,T=t\right)$. To establish the pointwise consistency of the estimator \eqref{estim}, we assume that the following general conditions hold:
\begin{enumerate}[label=(\textbf{A\arabic*})]
\item\label{exist_ass} There exists a $j$ such that  $|f_{U_j}\ft(v)|\neq 0$ for all $v\in \R$.
\item \label{KL_ass} i)The kernel function $L$ is such that $L\in L_1(\R)\cap L_2(\R)$, $L\ft$ is supported on $[-1,1]$; $|L\ft(v)-1|\downarrow0$ as $v \to 0$.\\
ii) $K\in L_1(\R)$ is bounded and is twice order kernel function.
\item\label{vanish_ass} $f_{(X,T)}(x,t)\neq 0$ for any $(x,t)$.
\item\label{bound_ass} $\tau^2f_{(X,T)}$, $f_{X}$, $f_T$, $rf_{(X,T)}$ and $r^2f_{(X,T)}$ are bounded and continuous. Moreover, $(mf_{(X,T)})\ft$, $f_{(X,T)}\ft\in L_1(\R^2)$. 
\item\label{var_ass} $\displaystyle\sum_{j=1}^n|f_{U_j}\ft(v)|^2\to \infty$, as $n\to \infty$ and $\displaystyle\inf_n\sum_{j=1}^n|f_{U_j}\ft(v)|^2>0$ $\forall v$.
\end{enumerate}

Assumptions~\ref{exist_ass}–\ref{bound_ass} are classical conditions in deconvolution and error-in-variables problems (see, for instance, \cite{FT93}). As pointed out by \cite{DM07}, Assumption~\ref{exist_ass} represents a key distinction between the homoscedastic and heteroscedastic frameworks in ensuring that the estimator \eqref{estim} is well defined. In the homoscedastic case, the Fourier transform of the error density, $f_{U}\ft$, is usually assumed to be nowhere vanishing. In contrast, under heteroscedasticity, it suffices to require that at least one of the Fourier transforms $f_{U_j}\ft$ does not vanish.
Assumption~\ref{var_ass} is a technical condition introduced to guarantee that the estimator \eqref{estim} is well defined. More precisely, it ensures that a sufficiently large number of the Fourier transforms $f_{U_j}\ft$ are nonzero, so that the denominator appearing in the definition of the deconvolution kernel does not degenerate.

%
The following theorem studies the pointwise consistency of the estimator \eqref{estim}. In order to establish its strong consistency, we further require the bandwidths $h$ and $b$ to satisfy the following more restrictive conditions.\\
\noindent The bandwidths $b$ and $h$ are such that $b=Ah$ for some $A>0$ and there exist $\delta>1$ and $\kappa>0$ such that 
\begin{equation}\label{band_ass}
c n^{(1-\delta+\kappa)/2}\leq b\to 0 \text{ and }\displaystyle\int_{|v|\leq 1/b} \left(\sum_{j=1}^n|f_{U_j}\ft(v)|^2\right)^{-2}dv=O(bn^{-\delta}).
\end{equation}
Note that if $b=Ah$ for some $A>0$ and for each $\eta>0$, $\inf_{|v|\leq \eta}|f_{U_j}\ft(v)|^2\geq Cn^{\alpha}$ with $\alpha>\delta/2$ and $C>0$, then \eqref{band_ass} is satisfied as soon as $\kappa=(\delta-1)/2$.\\
Let us now state the pointwise consistency of our estimator. 
\begin{theorem}\label{pointwise}
Assume that Assumptions \ref{exist_ass}--\ref{var_ass} holds. 
\begin{enumerate}
\item If $h\to0$, $b\to 0$ and $\int_{|v|\leq 1/b} \left(\sum_{j=1}^n|f_{U_j}\ft(v)|^2\right)^{-1}dv=o(h)$, then
\begin{eqnarray}
\label{weak}
\widehat{r}_n(x,t)\xrightarrow{\mathbb{P}}r(x,t) \text{  as } n \to \infty.
\end{eqnarray}
\item Moreover, assume that \eqref{band_ass} holds and $\|r\|_{\infty}<\infty$. If $\Eb|\varepsilon_j|^\ell\leq C_\ell<\infty$ for all integers $j$ and all $0<\ell\leq 2\lceil1/\kappa\rceil+2$, then
\begin{eqnarray}
\label{a.s.}
\widehat{r}_n(x,t)\xrightarrow{\text{a.s.} }r(x,t) \text{  as } n \to \infty.
\end{eqnarray}
where $\lceil a\rceil$ denotes the smallest integers larger or equal to $a$.
\end{enumerate}
\end{theorem}
Next, we want to establish the pointwise rates of convergence of the estimator $\rw(x,t)$ at a fixed arbitrarly $(x,t)\in \R^2$ and also show that these rates are optimal in a minimax sense with respect to any regression estimator in model \eqref{mod_het}. 
For any arbitrary fixed $(x,t)\in \R^2$, define the following class of functions:
\[
\begin{aligned}
\mathcal{F}_{\beta,C}
&=\Big\{ g \in \mathcal{C}^0(\R^2) \text{ s.t. $g$ and $g\ft$ are integrable and }
\left|
\frac{\partial^{k+\ell} g(x,t)}
     {\partial x^k \partial t^\ell}
\right|
\le C,\; 0 \le k+\ell = \beta, \\
&\qquad \forall~ (x,t)\in [ a, b]\times[c,d] \Big\},
\end{aligned}
\]
for some real numbers $a<b$ and $c< d$.
\begin{enumerate}[label=(\textbf{B\arabic*})]
\item\label{B1} $rf_{(X,T)}$ is continuous and integrable on $\R^2$, and $rf_{(X,T)}\in \mathcal{F}_{\beta,C_1}$ for some $C_1, \beta>0$.
\item\label{B2} $f_{(X,T)}\in \mathcal{F}_{\beta,C_1}$ with $\beta$ and $C_1$ defined as in \ref{B1}. 
\item\label{B3} $\|r^2f_{(X,T)}\|_{\infty}\leq C_2$, $\|f_{(X,T)}\|_{\infty}\leq C_3$ and $\|\tau^2f_{(X,T)}\|_\infty\leq C_4$ are bounded. 
\item\label{B4} Both $K$ and $L$ are $k$-th order kernel functions for some $k\geq 2\beta$, and $L\ft$ is supported on $[-1,1]$, with $\beta$ as in \ref{B1}. 
\item\label{B5} $\|r\|_\infty\leq C_5$ and $f_{(X,T)}(x,t)\geq C_6>0$.
\end{enumerate}
As pointed in \cite{DM07}, we also need some conditions on the error distribution. We suppose that there exist $\alpha$, $C>0$ and some monotonously decreasing functions $\overline{\varphi}_{j,n}(v)$ and $\underline{\varphi}_{j,n}(v)$ such that the following conditions hold:
\begin{enumerate}[label=(\textbf{C\arabic*})]
\item\label{C1}$\Pb\left(|U_j|\leq \alpha\right)\geq C$, $\forall j, n$.
\item\label{C2} $\big|f_{U_j}\ft(v)\big|\geq \underline{\varphi}_{j,n}(M)$, $\forall |v| \leq M$.
\item\label{C3} $\underline{\varphi}_{j,n}(v)\leq\big|f_{U_j}\ft(v)\big| \leq \overline{\varphi}_{j,n}(v)$, $\forall v>M$.
\item\label{C4} $\big|{f_{U_j}\ft}'(v)\big|\leq \overline{\varphi}_{j,n}(v)$, $\forall v >M$.
\item\label{C5}$\big|\underline{\varphi}_{j,n}(v)\big|\geq c_1\overline{\varphi}_{j,n}(c_2v)$, $\forall v>0$.
\end{enumerate}
Under these conditions, we assume that $M\geq 0$, $c_1>0$ and $c_2\geq 1$ are constants independent of both $j$ and $n$. As pointed out again by \cite{DM07}, condition \ref{C1} imposes a regularity requirement that prevents the error densities $f_{U_j}$ from becoming  excessively diffuse or overly smooth as $j$ increases. Moreover, conditions \ref{C1}--\ref{C5} constitute a weak form of monotonicity for $|f_{U_j}|$. \\
Finally, let $\mathcal{F}$ denote the class of all pairs $(r,f_{X,T})$ satisfying Conditions \ref{B1}--\ref{B5} and \ref{C1}--\ref{C5} with uniform constants and parameters. In what follows, the symbol \enquote{const} will denote an arbitrary positive constant.
\begin{theorem}\label{pointwise_rate}
Let $(x,t)\in\mathbb{R}^2$ be arbitrary fixed and conditions \ref{B1}--\ref{B5} and \ref{C1}--\ref{C5} hold. Assume that there exists a sequence $a_n\uparrow\infty$ such that  for some $C_8\geq C_7>0$, $\beta>1/2$,
\begin{eqnarray}\label{approx}
C_7a_n^{1+2\beta}\leq \sum_{j=1}^n\big|\overline{\varphi}_{j,n}(a_n)\big|^2\leq C_8a_n^{1+2\beta}
\end{eqnarray}
is valid for all $n$.
Then,
\begin{enumerate} 
\item When putting $h=cb$ with $c>0$ and $b=c_2a_n^{-1}$ with $c_2$ defined as in \ref{C5}, the estimator $\widehat{r}_n$ satisfies
\begin{eqnarray}
\label{limsup}
\limsup_{n\to \infty} \sup_{(r,f_{(X,T)})\in \mathcal{F}}\mathbb{P}\left[\left|\widehat{r}_n(x,t)-r(x,t)\right|^2>da_n^{1-2\beta}\right]\leq \const \cdot~ d^{-1} ~~\forall d>0.
\end{eqnarray}
\item For an arbitrary estimator $\widetilde r_n(x,t)$ and for sufficiently large constant $C$ in $\mathcal{F}=\mathcal{F}_{\beta,C}$, with $\beta>1/2$, there is some $C_{9}>0$ such that 
\begin{eqnarray}
\label{limsup}
\liminf_{n\to \infty} \sup_{(r,f_{(X,T)})\in \mathcal{F}}\mathbb{P}\left[\left|\widetilde{r}_n(x,t)-r(x,t)\right|^2>C_{9}a_n^{1-2\beta}\right]\geq \const.
\end{eqnarray}
\end{enumerate}

\end{theorem}

Theorem \ref{pointwise_rate} shows that the convergence rate of the proposed partial kernel deconvolution estimator \eqref{estim} attains the optimal convergence rate.

\section{Numerical studies}\label{numeric}

In order to evaluate the finite sample performance of our proposed kernel regression estimate, we consider in this section simulation studies in regression models with both supersmooth and ordinarly smooth measurements errors. Two different regression models are considered:
\begin{itemize}
\item Model 1:
\begin{eqnarray}
\label{modele1}
Y=r(X,T)+\varepsilon~~ \text{with}~~r(x,t)=x^2\exp{(-t^2/2)},
\end{eqnarray}
\item Model 2:
\begin{eqnarray}
\label{modele2}
Y=r(X,T)+\varepsilon~~\text{with}~~~ r(x,t)=x\theta+\cos(t)~~~\text{and}~~~\theta=3.
\end{eqnarray}
\end{itemize}
For each model, we took $X$ and $T$ as uniformly distributed on $[-2,2]$ while $\varepsilon\leadsto\mathcal{N}(0,0.25^2)$.
Our estimator was calculated  using the standard normal kernel function $K$ and another kernel $L$ whose characteristic function is 
$$
 L\ft(v)=(1-v^2)^3\mathds{1}_{\{|v|\leq 1\}}.
$$
For each model, we apply two different error models: $U_j\leadsto\text{Laplace}(\sigma_j)$ and $U_j\leadsto\mathcal{N}(0, \sigma_j^2)$, with $\var(U_j)=\sigma^2(1+j/n)$ for $j\in\{1, \ldots,n\}$, with $\sigma^2=0.2\var(T)$ and consider also the naive estimator by ignoring error using gaussian kernel for $K$ and $L$.  The sample sizes  considered is $n=100$, $500$ and $800$ and the number of simulations is $N=100$. For each simulation, we compute the average squared error (ASE) of our proposed estimator and the naive estimator (by ignoring error measurement) at the evaluated $50\times 50$ grid points in both directions from $[-2,2]$ to $[-2,2]$ using a bivariate $100$ bandwidths $(h,b)$ ranging from $[0.02,0.2]\times [0.02,0.2]$.The optimal bandwidths are selected to minimize the ASE among these $100$ pairwise candidates. \\
Noting that the reason of choosing a simple partial linear model \eqref{modele2} is the separability between the variables $x$  and $t$. When exploiting its structure, one can also estimate $r(x,t)$ by
\begin{eqnarray}
\label{rntilde}
\widetilde{r}_n(x,t)= x\widehat{\theta}_n+ \frac{\sum_{j=1}^{n}L_{U_j}\left(\frac{t-W_j}{b}\right) (Y_j-X_j\widehat{\theta}_n)}{\sum_{j=1}^{n}L_{U_j}\left(\frac{t-W_j}{b}\right)},
\end{eqnarray}
where $\widehat{\theta}_n$ is any $\sqrt{n}$- consistent estimator of $\theta$ and $L_{U_j}$ is defined in \eqref{lu_j}. For comparison in Model \ref{modele2}, both estimator $\rw_n(x,t)$ defined in \eqref{estim} and $\widetilde{r}_n(x,t)$ where considered in the simulation study. \\

\begin{table}[htb]
\begin{tabular}{lcccccl}
\toprule 
& \multicolumn{3}{c}{\textbf{Laplace error}} & \multicolumn{3}{c}{\textbf{Normal error}} \\
\cmidrule(lr){2-4}\cmidrule(lr){5-7}
& $\mathbf{n=100}$ & $\mathbf{n=500}$& $\mathbf{n=800}$& $\mathbf{n=100}$ & $\mathbf{n=500}$ & $\mathbf{n=800}$\\
\bottomrule
$\rw(x,t)$ &  0.0332   &0.0107& 0.0081&  0.0300&0.0081  &0.0053 \\
Naive &  0.0517&0.0111& 0.0091&0.0457  &0.0097 &0.0061 \\
\bottomrule
\end{tabular}
\caption{Average (over 100 simulation runs) squared error for estimating the model \ref{modele1}.}
\label{tab:ase_model1}
\end{table}

\begin{table}[htb]
\label{table1}
\begin{tabular}{lcccccl}
\toprule 
& \multicolumn{3}{c}{\textbf{Laplace error}} & \multicolumn{3}{c}{\textbf{Normal error}} \\
\cmidrule(lr){2-4}\cmidrule(lr){5-7}
& $\mathbf{n=100}$ & $\mathbf{n=500}$& $\mathbf{n=800}$& $\mathbf{n=100}$ & $\mathbf{n=500}$ & $\mathbf{n=800}$\\
\bottomrule
$\rw_n(x,t)$ &  0.0838   &0.0260&0.0129& 0.0780& 0.0192   & 0.0118\\
$\widetilde{r}_n(x,t)$  &0.0222& 0.0140& 0.0136&  0.0108& 0.0052& 0.0047\\
Naive & 0.1372 &0.0269& 0.0186&  0.1204 &0.0230& 0.0136\\
\bottomrule
\end{tabular}
\caption{Average (over 100 simulation runs) squared error for estimating the model \ref{modele2}.}
\label{tab:ase_model2}
\end{table}

As illustrated on Tables~\ref{tab:ase_model1} and \ref{tab:ase_model2}, the deconvolution method is robust to the error assumption and is significantly better than the naive estimate. According to the error distibutions, we can also see that the ASEs are relatively comparable.    \\

%

\begin{figure}[H]
\centering
\begin{subfigure}{0.48\textwidth}
    \centering
    \includegraphics[width=\linewidth]{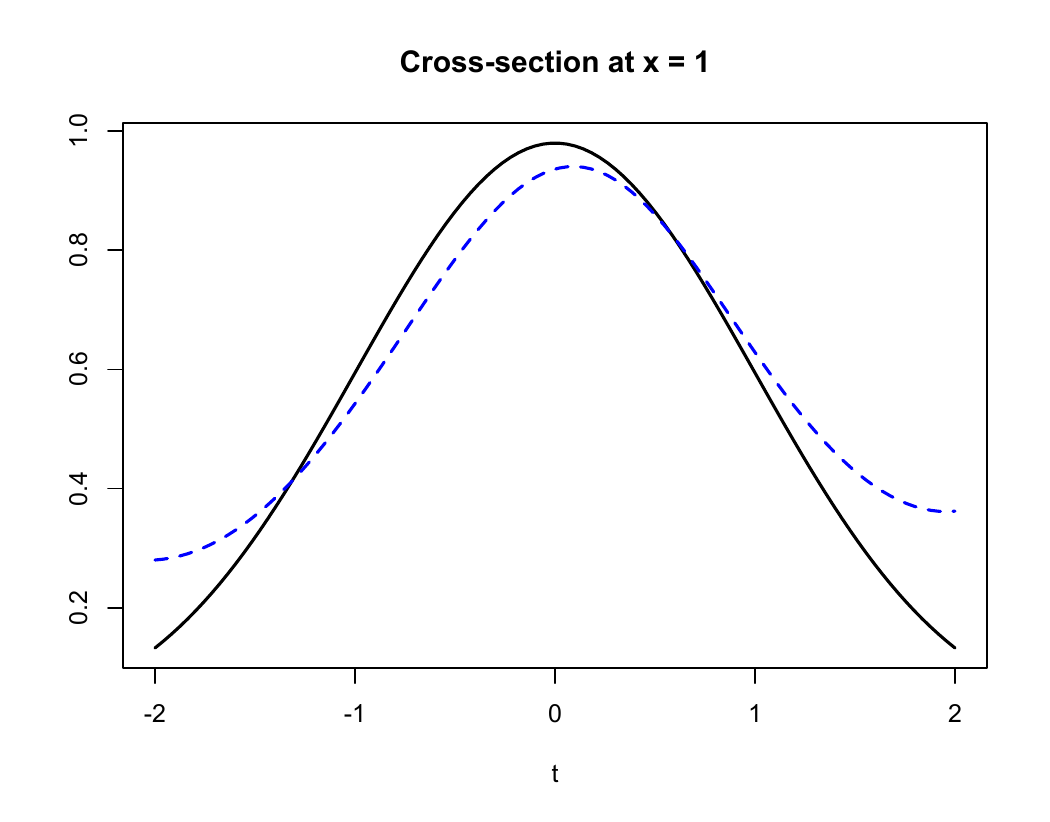}
\end{subfigure}
\hfill
\begin{subfigure}{0.46\textwidth}
    \centering
    \includegraphics[width=\linewidth]{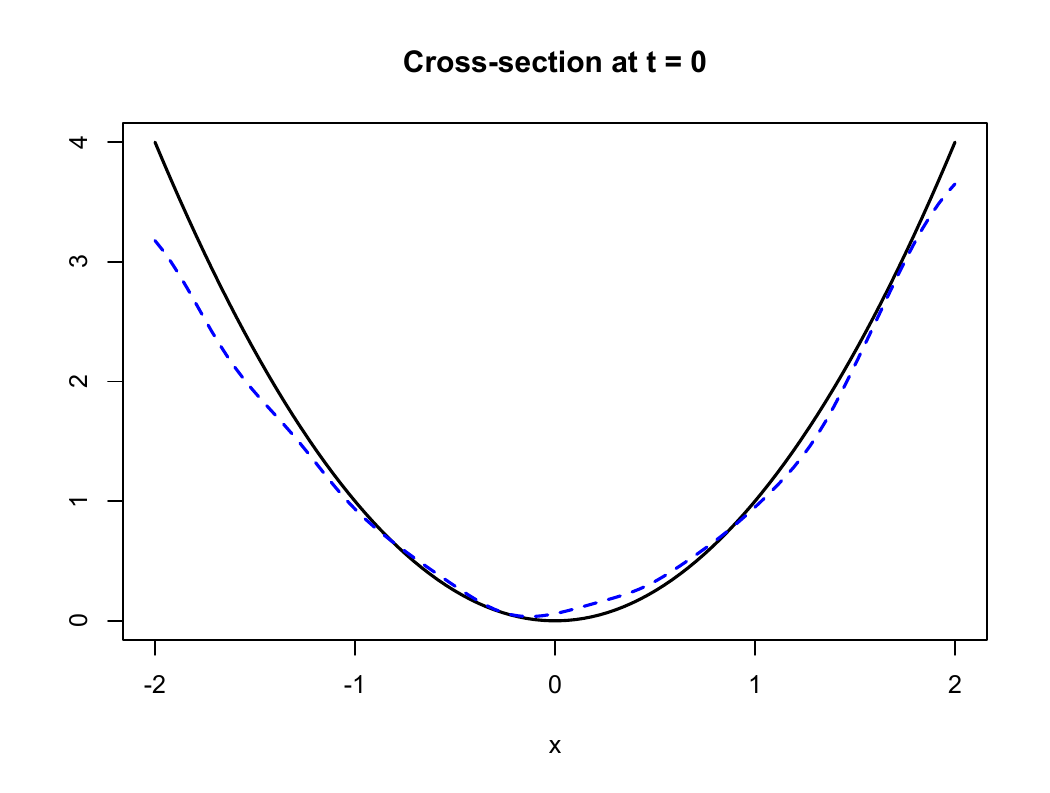}
\end{subfigure}
  \caption{Estimation of cross-sections at $x=1$ (on the left) and at $t=0$ (on the right)  from model \ref{modele1} under the Normal contaminated error.  The dashed line is the estimated curve $\rw_n$ and the solid line is the true curve. The sample size is  $n=500$. }
\label{fig:mod1}
\end{figure}

\begin{figure}[H]
\centering
\begin{subfigure}{0.48\textwidth}
    \centering
    \includegraphics[width=\linewidth]{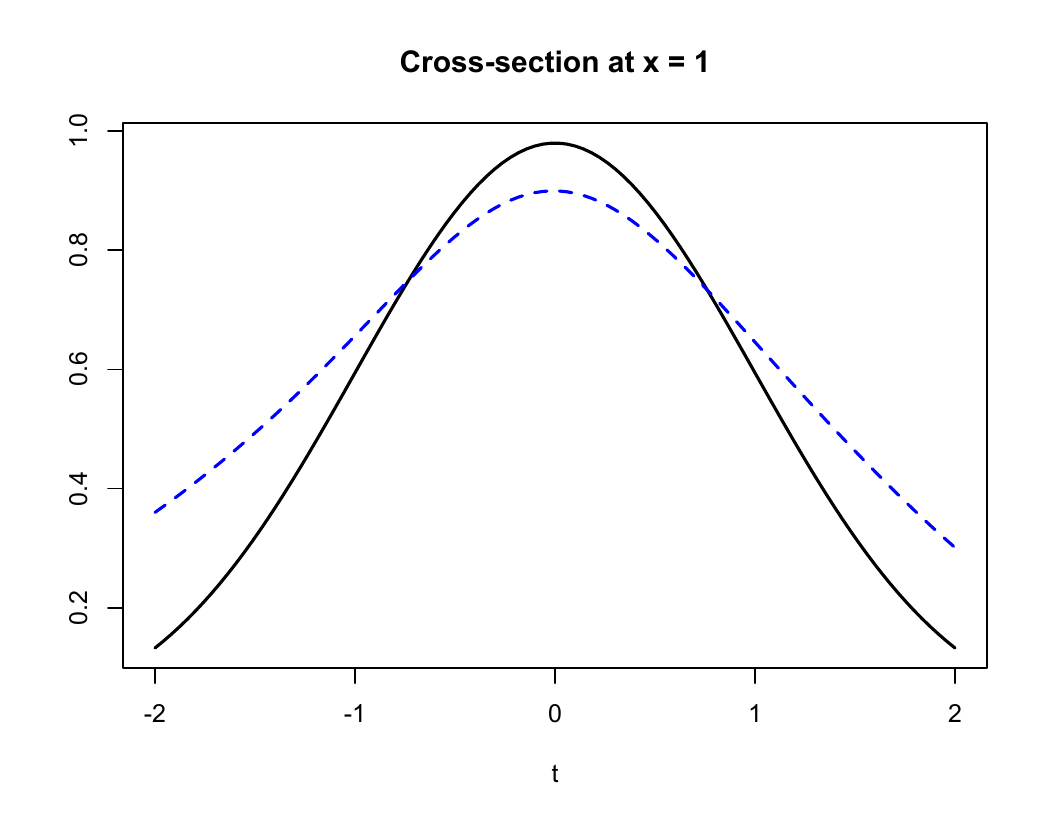}
\end{subfigure}
\hfill
\begin{subfigure}{0.46\textwidth}
    \centering
    \includegraphics[width=\linewidth]{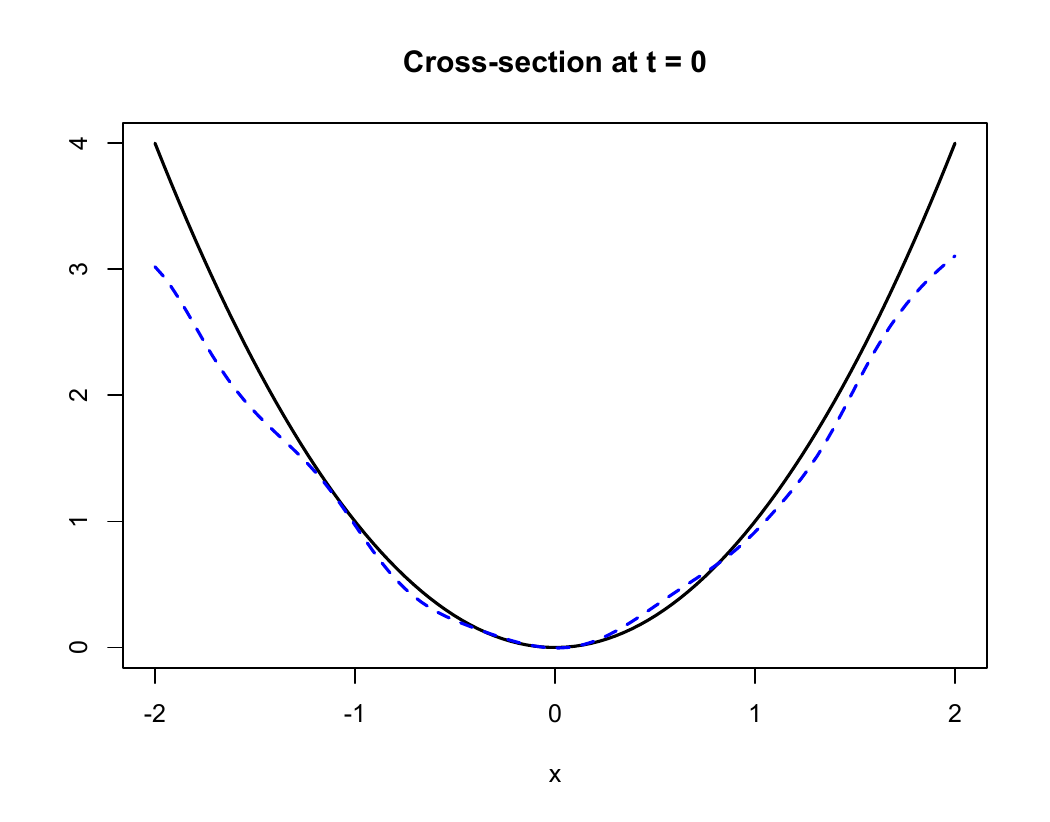}
\end{subfigure}
  \caption{Estimation of cross-sections at $x=1$ (on the left) and at $t=0$ (on the right)  from model \ref{modele1} under the Laplace contaminated error.  The dashed line is the estimated curve $\rw_n$ and the solid line is the true curve. The sample size is  $n=500$. }
\label{fig:mod1}
\end{figure}

\begin{figure}[H]
\centering
\begin{subfigure}{0.48\textwidth}
    \centering
    \includegraphics[width=\linewidth]{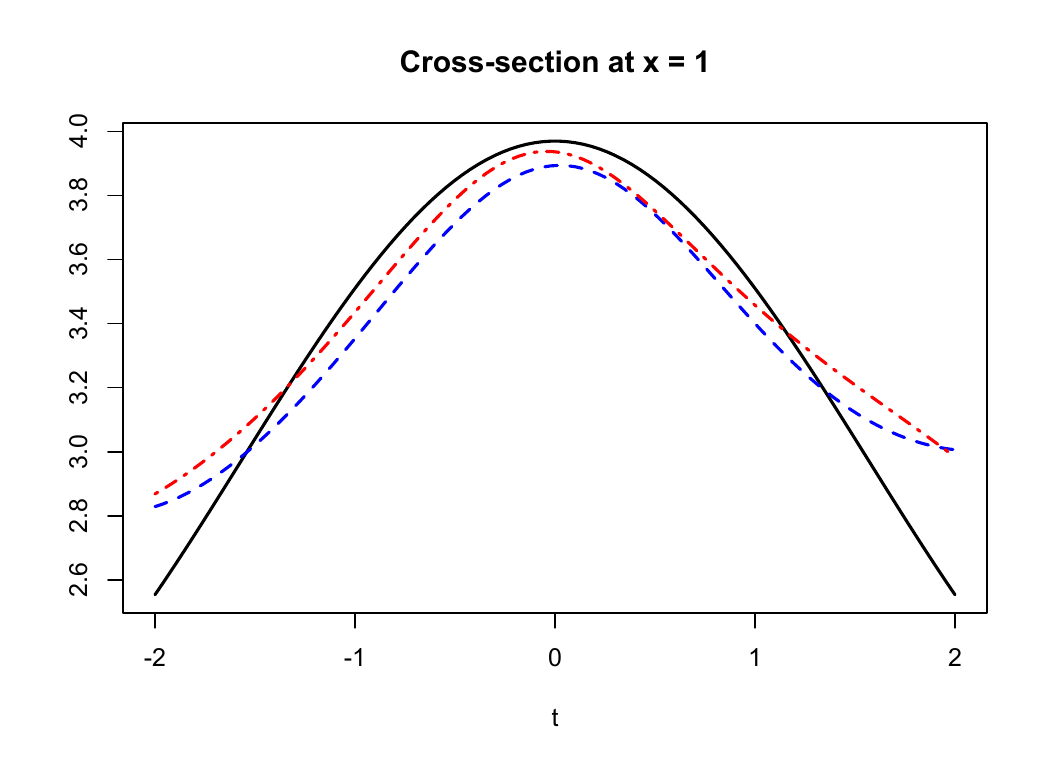}
\end{subfigure}
\hfill
\begin{subfigure}{0.48\textwidth}
    \centering
    \includegraphics[width=\linewidth]{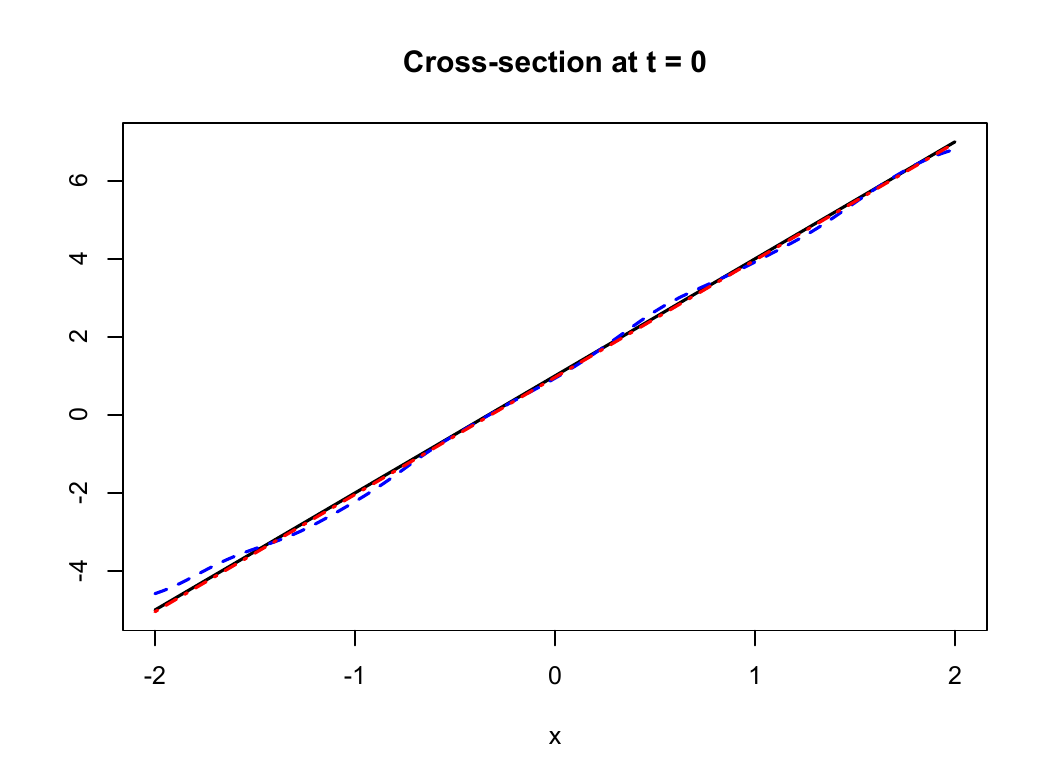}
\end{subfigure}
  \caption{Estimation of cross-sections at $x=1$ (on the left) and at $t=1$ (on the right)  from model \ref{modele2} under the Normal contaminated error.  The dashed line is the estimated curve $\rw_n$, the dotdashed line is $\widetilde r_n$ defined in \eqref{rntilde} and  the solid line is the true curve.The sample size is  $n=500$. }
\label{fig:mod2}
\end{figure}

\subsection*{Conclusion}
This paper studied nonparametric kernel regression estimation in a setting where some explanatory variables are measured with errors, and some are measured precisely. We consider a more general framework of the heteroscedastic errors measurement and shown that the proposed estimators achieve the optimal convergence rate under some regularity conditions.  Numerical experiments illustrated the finite-sample performance of the proposed estimation procedure. Bandwidth selection which is known to be  generally extremely difficult in the error-in-variables problem (even in the homoscedastic case) was not addressed in this work. The development of fully data-driven bandwidth selection methods is left for future research. Another direction for further investigation is the extension of the proposed methodology to settings where the measurement error density is unknown. 
\section{Proofs}\label{proofs}
Define
\begin{eqnarray}
\mw_n(x,t)=\sum_{j=1}^{n}Y_{j}K_h\left(x-X_j\right)L_{U_j}\left(\frac{t-W_j}{b}\right) \text{ and  }m(x,t)=r(x,t)f_{(X,T)}(x,t),
\end{eqnarray}
with $L_{U_j}$ defined in \eqref{lu_j}. Also, let 
$K_h(\cdot)=\frac{1}{h}K\left(\frac{\cdot}{h}\right)$ and $L_b(\cdot)=\frac{1}{b}L\left(\frac{\cdot}{b}\right)$.

\begin{lemma}\label{mse_num}
Assume that Assumptions \ref{exist_ass}--\ref{bound_ass} hold. Then
\begin{eqnarray*}
\Eb\left|\mw_n(x,t)-m(x,t)\right|^2=B_n^2+V_n,
\end{eqnarray*}
with
\begin{eqnarray}
\label{Bn2}
B_n^2 &=& \left|  \int_{\R^2} K_h \left(x-y\right) L_b\left(t-z\right)m(y,z)dydz-m(x,t)\right|^2=o((h^2+b^2)^2),\\
V_n 
&\leq &\left( \|\tau^2f_{X,T}\|_{\infty}+\|r^2f_{X,T}\|_{\infty}\right)\frac{1}{2\pi h}\int_\R\left| L\ft(vb)\right|^2 \left(\sum_{k=1}^n|f_{U_k}\ft(v)|^2\right)^{-1}dv.
\label{Vn}
\end{eqnarray}
\end{lemma}
\subsection*{Proof of Lemma \ref{mse_num}}
For any bivariate function $g(x,t)$, let $g_{(1)}\ft(v,t)=\int_\R \exp(iv x)g(x,t)dx$ be a fourrier transform of $g$ in $x$ only and $g_{(2)}\ft(x,v)=\int_\R \exp(iv t)g(x,t)dt$ be a fourrier transform of $g$ in $y$ only. 
For the bias term, by the independence between $U_j$ and $(X_j,T_j,Y_j)$, we have
\begin{eqnarray*}
\Eb(\widehat{m}_n(x,t))&= &(2\pi)^{-1}\sum_{j=1}^n\int_\R\exp(-ivt)\Eb\left[r(X_j,T_j)K_h\left(x-X_j\right)\exp(iv T_j)\right]\Eb\left[\exp(iv U_j)\right]L\ft(vb)\psi_j(v)dv\\
&=&
(2\pi)^{-1}\sum_{j=1}^n\int_\R\exp(-ivt)\Eb\left[r(X_j,T_j)K_h\left(x-X_j\right)\exp(iv T_j)\right]f_{U_j}\ft (v)L\ft(vb)\psi_j(v)dv.
\end{eqnarray*}
Note that 
\begin{eqnarray*}
\Eb\left[r(X_j,T_j)K_h\left(x-X_j\right)\exp(iv T_j)\right]
&= &  \int_{\R^2} K_h\left(x-\tilde x\right)m(\tilde x,\tilde t)\exp(iv \tilde t)d\tilde x d \tilde t\\
&=& \int _{\R} K_h\left(x-\tilde x\right)\left[\int_\R \exp(iv \tilde t)m(\tilde x,\tilde t)d \tilde t\right]d \tilde x\\
&:=&\int_\R K_h\left(x-\tilde x\right)m_{(2)}\ft(\tilde x,v)d\tilde x.
\end{eqnarray*}
We deduce that 
\begin{eqnarray*}
\Eb(\widehat{m}_n(x,t))&=&(2\pi)^{-1}\int_{\R^2}\exp(-ivt)K_h\left(x-\tilde x\right)m_{(2)}\ft(\tilde x,v)L\ft(vb)\sum_{j=1}^nf_{U_j}^{\text{th}}(v)\psi_j(v)dvd\tilde x\\
&=& \int_\R K_h\left(x-\tilde x\right)\left[(2\pi)^{-1}\int_\R\exp(-ivt)m_{(2)}\ft(\tilde x,v)L\ft(vb)dv\right]d\tilde x.
\end{eqnarray*}
Since $L\ft(vb)=L_b\ft(v)$ and $(L_b*m(\cdot,v))\ft=m_{(2)}\ft(\cdot,v) L_b\ft(v)$, then by Fourier inversion, we have
$$
(2\pi)^{-1}\int_\R\exp(-ivt)m_{(2)}\ft(\tilde x,v)L\ft(vb)dv=(2\pi)^{-1}\int_\R\exp(-ivt)(L_b*m(\cdot,v))\ft dv=L_b*m(\cdot,t).
$$
It follows that 
\begin{eqnarray*}
\Eb(\widehat{m}_n(x,t))&=& \int_\R K_h\left(x-\tilde x\right)L_b*m(\tilde x, t)d\tilde x\\&=&\int_{\R^2} K_h\left(x-\tilde x\right)L_b\left(t-\tilde t\right)m(\tilde x,\tilde t)d\tilde x d\tilde t=\int_{\R^2}K(u)L(v)m(x-uh,t-vb)dudv.
\end{eqnarray*}
Therefore, Taylor's expansion combined with assumption \ref{KL_ass} imply \eqref{Bn2} in Lemma \ref{mse_num}. For the variance term, note that $\left|\Eb\left[Y_j^2\vert X_j=x,T_j=t\right]\right|\leq \|\tau^2f_{X,T}\|_{\infty}+\|r^2f_{X,T}\|_{\infty}:=C$.
Using the independence between $Y_j|(X_j,T_j)$ and $U_j$ on one hand, and $X_j$ and $T_j$ on the other hand, we have
\begin{eqnarray*}
V_n&=&\sum_{j=1}^n\var\left(Y_jK_h(x-X_j)(2\pi)^{-1}\int_{\R}\exp(-ivt)\exp(ivW_j)L\ft(vb)\psi_j(v)dv\right)\\
&\leq & \sum_{j=1}^n\Eb\left\{\Eb\left[Y_j^2\vert X_j,T_j\right]K_h^2(x-X_j)\left|(2\pi)^{-1}\int_{\R}\exp(-ivt)\exp(ivW_j)L\ft(vb)\psi_j(v)dv\right|^2\right\}\\
& \leq & C\sum_{j=1}^n\Eb\left[K_h^2(x-X_j)\right]\Eb\left[\left|(2\pi)^{-1}\int_{\R}\exp(-ivt)\exp(ivW_j)L\ft(vb)\psi_j(v)dv\right|^2\right].
\end{eqnarray*}
Next, since $W_j=T_j+U_j$, then $f_T*f_{U_j}(t)=\int_\R f_T(z)f_{U_j}(t-z)dz$ is the density function of $W_j$. Then, we have
\begin{eqnarray*}
\lefteqn{\Eb\left[\left|(2\pi)^{-1}\int_{\R}\exp(-ivt)\exp(ivW_j)L\ft(vb)\psi_j(v)dv\right|^2\right]}\\
 & &=\int_\R\left| (2\pi)^{-1} \int_\R \exp(-iv(t-z))L\ft(vb)\psi_j(v)dv\right|^2f_T*f_{U_j}(z)dz\\
 &=&\int_\R\left| (2\pi)^{-1} \int_\R \exp(-ivu)L\ft(vb)\psi_j(v)dv\right|^2f_T*f_{U_j}(t-u)du\\
&\leq &\|f_T*f_{U_j}\|_{\infty}\int_\R\left| (2\pi)^{-1} \int_\R \exp(-ivu)L\ft(vb)\psi_j(v)dv\right|^2du.
\end{eqnarray*}
Therefore, using Parseval's identity, 
 we can write
\begin{eqnarray*}
\Eb\left[\left|(2\pi)^{-1}\int_{\R}\exp(-ivt)\exp(ivW_j)L\ft(vb)\psi_j(v)dv\right|^2\right] &\leq &\|f_T*f_{U_j}\|_{\infty}(2\pi)^{-1}\int_\R |L\ft(vb)|^2|\psi_j(v)|^2dv.
\end{eqnarray*}
Noting that by \ref{bound_ass},  $f_T*f_{U_j}$ is bounded, it follows that 
\begin{eqnarray*}
V_n&\leq&C h^{-1}\int_\R K^2(z)f_X(x-hz)dz(2\pi)^{-1}\int_\R |L\ft(vb)|^2\sum_{j=1}^n|\psi_j(v)|^2dv  \\
&\leq & 
Ch^{-1}\int_\R K^2(z)f_X(x-hz)dz(2\pi)^{-1}\int_\R |L\ft(vb)|^2\left[\sum_{k=1}^n|f_{U_k}^{\text{th}}(v)|^2\right]^{-1}dv,
\end{eqnarray*}
which proves \eqref{Vn}. Observe that Lemma \ref{mse_num} focuses on estimation of $m(x,t)$. Setting $Y_j=1$ almost surely, we get immediately the same result for the estimation of $f_{(X,T)}(x,t)$ by $\widehat{f}_n(x,t)$ defined in \eqref{denom}.
\hfill$\square$
\begin{lemma}\label{mom_num}
Assume that Assumptions \ref{exist_ass}--\ref{var_ass} hold and also $\|r\|_{\infty}<\infty$. Then, for any integer $\gamma>0$, if $\Eb|\varepsilon_j|^\ell\leq C_\ell<\infty$ $\forall \ell\leq 2\gamma$, we have
\begin{eqnarray*}
\Eb\left|\widehat{m}_n(x,t)-m(x,t)\right|^{2\gamma}=O\left(h^{-\gamma}b^{-\gamma}n^\gamma\right)\left[h^{-1}\int_{|v|\leq 1/b} \left(\sum_{k=1}^n|f_{U_k}\ft(v)|^2\right)^{-2}dv\right]^\gamma.
\end{eqnarray*} 
\end{lemma}

\subsection*{Proof of Lemma \ref{mom_num}}
Let $Z_j(x,v)=Y_j\exp\left(ivW_j\right)K_h\left(x-X_j\right)-\Eb\left[Y_j\exp\left(ivW_j\right)K_h\left(x-X_j\right)\right]$ and note that by the independence between $(X_j,T_j,Y_j)$ and $U_j$, we can write 
\begin{eqnarray*}
\Eb\left|Z_j(x,v)\right|^{2\gamma}&\leq &O(1) h^{-2\gamma}\Eb\left[\left|Y_j\exp(ivT_j)\exp(ivU_j)K\left(\frac{x-X_j}{h}\right)\right|^{2\gamma}\right]\\ & \leq & O(1) h^{-2\gamma}\Eb\left[\left|Y_j\exp(ivT_j) K\left(\frac{x-X_j}{h}\right)\right|^{2\gamma} \right]\\&\leq &O(1) h^{-2\gamma}\Eb\left[\Eb\left[\left|Y_j\right|^{2\gamma}\vert X_j,T_j \right] \left|K\left(\frac{x-X_j}{h}\right)\right|^{2\gamma}\right]\\ & \leq &O(1) h^{1-2\gamma}\int_\R |K(z)|^{2\gamma}f_X(x-hz)dz=O(h^{1-2\gamma})=O(h^{-2\gamma}).
\end{eqnarray*}
 Then, we have
\begin{eqnarray*}
\Eb\left|\mw(x,t)-m(x,t)\right|^{2\gamma}&=&\Eb\left|\sum_{j=1}^n(2\pi)^{-1}\int_\R \exp(-ivt)L\ft(vb)Z_j(x,v)\psi_j(v)dv\right|^{2\gamma}\\ & \leq & O(1)\sum_{j_1}\ldots\sum_{j_{2\gamma}}\int_\R\ldots\int_\R\prod_{k_1=1}^\gamma |L\ft(v_{k_1}b)||\psi_{j_{k_1}}(v_{k_1})|
\\ & & \times \prod_{k_2=\gamma+1}^{2\gamma}|L\ft(v_{k_2}b)||\psi_{j_{k_2}}(v_{k_2})|\mathds{1}_{\{\#\{j_1, \ldots, j_{2\gamma}\}\leq \gamma\}}\left[\sup_{k}\Eb\left|Z_k(x,v)\right|^{2\gamma}\right]dv_1\ldots dv_{2\gamma}\\
&\leq&O(1) h^{-2\gamma} \sum_{j_1}\ldots\sum_{j_{2\gamma}} \mathds{1}_{\{\#\{j_1, \ldots, j_{2\gamma}\}\leq \gamma\}}\prod_{\ell=1}^{2\gamma}\int_{\R} |L\ft(vb)||\psi_{j_\ell}(v)|dv\\
& \leq & O(1) h^{-2\gamma} \sum_{j_1}\ldots\sum_{j_{2\gamma}} \mathds{1}_{\{\#\{j_1, \ldots, j_{2\gamma}\}\leq \gamma\}}\left[\int_{|v|\leq 1/b} \left(\sum_{k=1}^n|f_{U_k}\ft(v)|^2\right)^{-2}dv\right]^\gamma \\
& \leq& O(1) n^\gamma h^{-2\gamma}b^{-\gamma} \left[\int_{|v|\leq 1/b} \left(\sum_{k=1}^n|f_{U_k}\ft(v)|^2\right)^{-2}dv\right]^\gamma,
\end{eqnarray*}
where we  used the fact that $\Eb(Z_{j_1}\cdots Z_{j_{2\gamma}})=0$ contains more than $\gamma$ different elements because $Z_j$ is centered and also we applied the Cauchy-Schwarz inequality to the integrals. \hfill$\square$
\subsection*{Proof of Theorem \ref{pointwise}}
For the proof of the first part of Theorem \ref{pointwise},  by Assumptions \ref{KL_ass} and \ref{bound_ass}, we have  $B_n^2$ tends to 0 as $(h,b)\to(0,0)$. Since $|L\ft|$ is bounded, conditions of Theorem \ref{pointwise} combined with Lemma \ref{mse_num} ensure that the variance term $V_n$ also converges to. Therefore, we have the convergence  in probability of $\mw_n(x,t)$ to $m(x,t)$ and $\fw_n(x,t)$ to $f_{(X,T)}(x,t)$. By usual arguments and Assumption \ref{vanish_ass}, we deduce that $\rw_n(x,t)$ has the weak limit $m(x,t)/f_{(X,T)}(x,t)$. \\
For the second part of Theorem \ref{pointwise}, we consider
\begin{eqnarray}
\label{decomp}
\left|\mw_n(x,t)-m(x,t)\right|\leq \left|\mw_n(x,t)-\Eb(\mw_n(x,t))\right|+ \left|\Eb(\mw_n(x,t))-m(x,t)\right|.
\end{eqnarray}
Note that, as shown in the first part, the deterministic term $B_n$  converges to 0 as $(h,b)\to (0,0)$. For the first term in \eqref{decomp}, we have, for any $\varepsilon>0$,
\begin{eqnarray*}
\sum_{n=1}^\infty \Pb\left[\left|\mw_n(x,t)-\Eb(\mw_n(x,t))\right|>\varepsilon\right]&\leq &\varepsilon^{-2\gamma}\sum_{n=1}^\infty \Eb\left|\mw_n(x,t)-\Eb(\mw_n(x,t))\right|^{2\gamma}\\ 
&\leq &O(1)\varepsilon^{-2\gamma}\sum_{n=1}^\infty O(h^{-\gamma} b^{-\gamma})n^\gamma\left[h^{-1}\int_{|v|\leq 1/b} \left(\sum_{k=1}^n|f_{U_k}\ft(v)|^2\right)^{-2}dv\right]^\gamma\\
&\leq &O(1)\varepsilon^{-2\gamma}\sum_{n=1}^\infty n^{-(1-\delta+\kappa)\gamma}n^{-\gamma(\delta-1)}=O(1)\varepsilon^{-2\gamma}\sum_{n=1}^\infty n^{-\gamma\kappa}<\infty,
\end{eqnarray*}
if we take $\gamma\geq \lceil1/\kappa\rceil+1$, we used Markov's inequality, Lemma \ref{mom_num} and \eqref{band_ass}.
Therefore, Borel-Cantelli Lemma implies that $|\mw_n(x,t)-m(x,t)|\xrightarrow{\text{a.s.} }0 \text{  as } n \to \infty.$ Similar arguments provide $|\fw_n(x,t)-f_{(X,T)}(x,t)|\xrightarrow{\text{a.s.} }0 \text{  as } n \to \infty$, which prove Theorem \ref{pointwise}.\hfill$\square$
\subsection*{Proof of Theorem \ref{pointwise_rate}}
Note that 
\begin{align*}
\rw_n(x,t)-r(x,t)
&= \frac{\mw_n(x,t)}{\fw_n(x,t)}-\frac{m(x,t)}{f_{(X,T)}(x,t)}
= \frac{\mw_n(x,t)-m(x,t)}{\fw_n(x,t)} - m(x,t)\frac{\fw_n(x,t)-f_{(X,T)}(x,t)}{f_{(X,T)}(x,t) \fw_n(x,t)}.
\end{align*}
Then, on the set $\{|\fw_n(x,t)-f_{(X,T)}(x,t)|\le C_6/2\}$, condition \ref{B5} implies that 
\[
|\rw_n(x,t)-r(x,t)| \le \frac{2}{C_6}|\mw_n(x,t)-m(x,t)| + \frac{2C_5}{C_6}|\fw_n(x,t)-f_{(X,T)}(x,t)|.
\]
Therefore,  
\begin{align*}
|\rw_n(x,t)-r(x,t)|^2
&\le \frac{8}{C_6^2}|\mw_n(x,t)-m(x,t)|^2
+ \frac{8C_5^2}{C_6^2}|\fw_n(x,t)-f_{(X,T)}(x,t)|^2.
\end{align*}
 Markov's inequality ensures that 
\begin{eqnarray*}
\lefteqn{\Pb\left(\left|\rw_n(x,t)-r(x,t)\right|^2>d a_n^{1-2\beta}\right)}\\
 & \leq & \Pb\left(\left|\mw_n(x,t)-m(x,t)\right|^2>(d C_6^2/16)a_n^{1-2\beta}\right)+ \Pb\left(\left|\fw_n(x,t)-f_{(X,T)}(x,t)\right|^2>d(c_6^2/16C_5^2)a_n^{1-2\beta}\right)\\ &&+\Pb\left(\left|\fw_n(x,t)-f_{(X,T)}(x,t)\right|^2>C_6^2/4\right)\\ 
 & \leq & \text{const. }d^{-1}a_n^{1-2\beta}\left\{ \Eb\left|\mw_n(x,t)-m(x,t)\right|^2+\Eb\left|\fw_n(x,t)-f_{(X,T)}(x,t)\right|^2 \right\}.
\end{eqnarray*}
Therefore, the first part of Theorem \ref{pointwise_rate} follows from \eqref{Vn} of Lemma \ref{mse_num}, \ref{C2}, \ref{C3}, \ref{C5} and \eqref{approx}.\\
\noindent For second part of Theorem \ref{pointwise_rate}, without loss of generality, let us estimate $r(x,t)$ at $x=t=0$.  We introduce the specific densities $f(t)=\left(1-\cos(t)\right)/(\pi t^2)$ with the Fourrier transform $f\ft(v)=(1-|v|)\mathds{1}_{[-1,1]}(v)$, the supersmooth Cauchy density $s_T(t)=1/\pi(1+t^2)$, $s_Y(y)=(1/2)\exp(-|y|)$, $\Delta_Y(y)=(1/4)\sign(y)\exp(-|y|)$.  Let $s_X(x)$ be the standard normal probability density function and $\Delta_T(t)=a_n^{-\beta}\cos(2a_nt)f(a_nt)$. As competing densities for $(X_j,T_j,Y_j)$, we give
$$
f_{(X,T,Y),p}(x,t,y)=s_X(x)s_T(t)s_Y(y)+p\const s_X(x)\Delta_T(t)\Delta_Y(y),
$$  
with $p\in \{0,1\}$ and $a_n$ as in \eqref{approx}. By using a symmetry argument for $\Delta(y)$, we obtain
 $$f_{(X,T),p}(x,t)=\int f_{(X,T,Y),p}(x,t,y)dy=s_X(x)s_T(t),$$ 
 and also 
$$
r_p(x,t)f_{(X,T)}(x,t)=\int yf_{(X,T,Y)}(x,t,y)dy=p\const a_n^{-\beta} f(a_nt)\cos(2a_nt)s_X(x). 
$$ 
Then, we can verify that $rf_{(X,T)}$ and $f_{(X,T)}$ fulfilled the smoothness conditions (B1) and (B2).
Let 
$$
h_{j,p}(x,w,y)=\int f_{(X,T,Y),p}(x,t,y)f_{U_j}(w-y)dy
$$
denote the density of $(X,T,Y)$.  Setting $D=\inf_n \left\{a_n^{2\beta-1}\left|r_0(0,0)-r_1(0,0)\right|^2/4\right\}$, then for $n$ large enough, the events $\left\{ \left|\widetilde r_n(0,0)-r_0(0,0)\right|^2< Da_n^{1-2\beta}\right\}$ and $\left\{ \left|\widetilde r_n(0,0)-r_1(0,0)\right|^2< Da_n^{1-2\beta}\right\}$ are disjoint. Hence, for any estimator $\tilde r_n$ of the regression function $r$, we have
\begin{eqnarray}
\label{min}
\lefteqn{2\sup_{(r,f_{(X,T)})}\Pb\left[\left|\widetilde r_n(0,0)-r(0,0)\right|^2\geq Da_n^{1-2\beta}\right]}\nonumber\\
& \geq & \Pb_{(r_0,f_{(X,T),0})}\left[\left|\widetilde r_n(0,0)-r_0(0,0)\right|^2\geq Da_n^{1-2\beta}\right]  +\Pb_{(r_1,f_{(X,T),1})}\left[\left|\widetilde r_n(0,0)-r_1(0,0)\right|^2\geq Da_n^{1-2\beta}\right] \nonumber\\& \leq&  \int \ldots \int \min\left(\prod_{j=1}^nh_{j,0}(x_j,w_j,y_j),\prod_{j=1}^nh_{j,1}(x_j,w_j,y_j)\right) dx_1dw_1dy_1\ldots dx_ndw_ndy_n.
\end{eqnarray}
We have to show that \eqref{min} is bounded away from 0. Using LeCam's inequality, this corresponds to 
\begin{eqnarray}
\label{prod}
\prod_{j=1}^n\int\int\int \sqrt{h_{j,0}(x,w,y) h_{j,1}(x,w,y)}dxdwdy\geq \const>0,
\end{eqnarray}
which is equivalent for $n$ large enough 
\begin{eqnarray*}
\sum_{j=1}^n\left|\ln \int\int\int \sqrt{h_{j,0}(x,w,y) h_{j,1}(x,w,y)}dxdwdy\right|\leq \const.
\end{eqnarray*}
Note that, from the definition of $f_{(X,T,Y),p}$, we have $h_{j,p}(x,w,y) \geq \const \cdot\exp{(-|y|)}g_j(w)s_X(x)$, where 
$g_j(w)=\int s_T(t)f_{U_j}(w-t)dt$ is a probability density function. Then, we get a positive lower bound on all integrals occurring in \eqref{prod} 

Using the inequality $|\ln(x)|\leq |(x-1)/x|$ for $x\in]0,1]$, we obtain
\begin{eqnarray*}
\lefteqn{
\sum_{j=1}^n\left|\ln \int\int\int \sqrt{h_{j,0}(x,w,y) h_{j,1}(x,w,y)}dxdwdy\right|}\\& &
\leq 
\const\cdot \sum_{j=1}^n\left(1-\int\int\int \sqrt{h_{j,0}(x,w,y) h_{j,1}(x,w,y)}dxdwdy \right)\\&&\leq \const\cdot\sum_{j=1}^n \chi^2\left(h_{j,0},h_{j,1}\right),
\end{eqnarray*}
where $\chi^2(f,g)=\int (f-g)^2/f dx$ denote the chi-squared distance and we use the fact that $1-\int \sqrt{fg}\leq \frac{1}{2}\chi^2(f,g)$ for any probability densities functions $f$ and $g$. Finally, \eqref{prod} follows if we show that 
\begin{eqnarray}
\label{chisq}
\chi^2\left(h_{j,0},h_{j,1}\right)=O(1).
\end{eqnarray}
For this, we have, for $p\in\{0,1\}$,
\begin{eqnarray*}
h_{j,p}(x,w,y)&=&s_X(x)s_Y(y)\int s_T(t)f_{U_j}(w-t)dt.\\ &\geq &s_X(x)s_Y(y)\int_{|t|\leq d}s_T(w-t)f_{U_j}(t)dt\\& \geq&s_X(x)s_Y(y) \frac{1}{2(1+2w^2+2d^2)}\int_{|t|\leq d}f_{U_j}(t)dt\\&\geq& c\exp(-|y|-x^2/2)\frac{1}{1+w^2},
\end{eqnarray*}
where condition \ref{C1} is used and selecting $d$ sufficiently large and $c$ sufficiently small. 
Indeed, \eqref{chisq} can be proved following the same lines as in the proof of (A.5) in \cite{DM07}, which conclude the proof of Theorem \ref{pointwise_rate}. \hfill$\square$

\end{document}